\documentclass[12pt,a4paper]{article}
\usepackage[utf8]{inputenc}
\usepackage{amsmath}
\usepackage{amsfonts}
\usepackage{color}
\usepackage{amssymb}
\usepackage{mathrsfs}
\usepackage{esint}
\usepackage[colorinlistoftodos]{todonotes}
\usepackage[colorlinks=true]{hyperref}
\hypersetup{
	colorlinks=blue,%
	citecolor=red,%
	filecolor=black,%
	linkcolor=blue,%
	urlcolor=blue
}
\usepackage[margin=1.3in]{geometry}

\usepackage[amsmath,thmmarks,hyperref]{ntheorem}
{
	\theoremstyle{nonumberplain}
	\theoremheaderfont{\bfseries}
	\theorembodyfont{\normalfont}
	\theoremsymbol{\mbox{$\Box$}}
	\newtheorem{pf}{Proof.}
}

\numberwithin{equation}{section}
\def\R{\mathbb{R}}
\def\S{\mathbb{S}}
\def\N{\mathbb{N}}
\def\H{\mathbb{H}}
\def \no{\nonumber}

\newcommand{\ud}{\mathrm{d}}
\def\pa{\partial}
\newtheorem{thm}{Theorem}[section]

\newtheorem{ex}{Example}[section]
\newtheorem{lem}{Lemma}[section]
\newtheorem{rem}{Remark}[section]

\newtheorem*{thm A}{Theorem A}
\newtheorem*{thm B}{Theorem B}
\newtheorem{cor}{Corollary}[section]
\usepackage{upgreek}
\usepackage{graphicx}

\allowdisplaybreaks
\title{article}
\begin{document}
	\title{\bf On some  rigidity  theorems of Q-curvature} 
	\date{\today}

\author{Yiyan Xu\thanks{Y. Xu: xuyiyan@nju.edu.cn;  $^\dag$S. Zhang: dg21210019@smail.nju.edu.cn}~ and Shihong Zhang$^\dag$
\footnote{Corresponding author: S. Zhang} }

\footnotetext[1]{The first author has been  supported in part by NSFC under Grant
	No. 11501285 and No. 11871265.}
	\maketitle

	{\noindent\small{\bf Abstract:}
In this paper, we investigate the rigidity of Q-curvature. Specifically, we consider a closed, oriented $n$-dimensional ($n\geq6$) Riemannian manifold $(M,g)$ and prove the following results under the condition $\int_{M}\nabla R\cdot\nabla \mathrm{Q}\ud V_g\leq0$.
(1) If $(M,g)$ is locally conformally flat with nonnegative Ricci curvature, then $(M,g)$ is isometric to a quotient of $\R^n$, $\S^n$, or $\R\times\S^{n-1}$.
(2) If $(M,g)$ has $\delta^2 W=0$ with nonnegative sectional curvature, then $(M,g)$ is isometric to a quotient of the product of Einstein manifolds.
Additionally, we investigate some rigidity theorems involving Q-curvature about hypersurfaces in simply-connected space forms. We also show the uniqueness of metrics with constant scalar curvature and constant Q-curvature in a fixed conformal class.

	\medskip 

{{\bf $\mathbf{2020}$ MSC:} 53C24, 53C18 (58C35, 58D17)}


\section{Introduction}

The Q-curvature, introduced by Tom Branson \cite{Branson85}, is defined in terms of the Ricci tensor ${\rm Ric}$ and scalar curvature $R$ of the manifold $(M,g)$ as follows:
\begin{equation}\label{Q-curvature}
\mathrm{Q}=-\frac{1}{2(n-1)}\Delta R-\frac{2}{(n-2)^2}| \mathrm{Ric}|^2+\frac{n^3-4n^2+16n-16}{8(n-1)^2(n-2)^2}R^2.
\end{equation}

In recent years,  the Q-curvature has attracted lots of attention since it arises naturally in conformal geometry in the context of conformally covariant operators--  the Paneitz operator (see \eqref{Paneitz}). We signal to the interested reader the survey  \cite{A. Chang and M.Eastwood, F.Hang and P. Yang 3} and also the references therein as a  review of the current literature, open problems, and recent developments. Moreover, an Einstein metric has constant Q-curvature. To some extent, a metric with constant Q is a generalization of Einstein metric.
In this note,   we will be mainly interested in the  rigidity property of Q-curvature with some additional conditions and we also note that Lin and Yuan have
found many interesting rigidity results in \cite{Lin Wei 1, Lin Wei 2, Lin Wei 3}.

It is worth mentioning that there has been extensive research on complete locally conformally flat manifolds with nonnegative Ricci curvature. 
Zhu \cite{S.Zhu}  proved  that the universal
covering of those manifolds either  is isometric to $\mathbb{R}\times \S^{n-1}$, or
is conformally equivalent to $\S^n$ or $\mathbb{R}^n$. Carro and Herzlich \cite{Carro&Herzlich} extended these results  and  achieved more   information on their  topology  and   geometry. More precisely, they  proved that 
those manifolds are either flat, or locally isometric to
 to $\mathbb{R}\times \S^{n-1}$,   or are globally conformally equivalent to $\mathbb{R}^n$ or to a
spherical space form $\S^n/\Gamma$.
If one adds extra assumptions such as constant  scalar curvature, Cheng \cite{Q. Cheng} proved some rigidity theorems.  This is one  motivation
for our paper to study the rigidity problem.  In particular,  Cheng proved the following theorem:

\begin{thm A}[\protect{Cheng \cite[Theorem 1]{Q. Cheng}}]
Let $(M^n, g)$ be a compact, oriented, locally conformally flat Riemannian manifold with constant scalar curvature. If  $\mathrm{Ric}\geq 0$, then $(M^n, g)$ is isometric to a quotient of either flat Euclidean space $\mathbb{R}^n$, round sphere $\mathbb{S}^n$, or cylinder $\mathbb{R}\times \mathbb{S}^{n-1}$.
\end{thm A}

In the following, we aim to establish a rigidity result for closed, oriented, locally conformally flat manifolds with non-negative Ricci curvature. To achieve this, we will introduce an integral condition involving the Q-curvature. It is important to emphasize that the previous result, as presented by \cite{Carro&Herzlich,Q. Cheng,S.Zhu}, did not specifically address the Q-curvature. Hence, our contribution extends Cheng's previous result by incorporating the Q-curvature into the analysis.
\begin{thm}\label{Thm 1.1}
Suppose that $ (M^n,g)$ is a closed, oriented, locally conformally flat $n$-dimensional Riemannian manifold with $n\geq6$ and the following conditions holds:
	$$\mathrm{Ric}\geq0,\quad\quad\quad \mathrm{and}\quad\quad\int_{M}\nabla R\cdot\nabla \mathrm{Q}\ud V_g\leq0.$$ Then, $(M^n,g)$ is isometric to a quotient of flat Euclidean space  $\R^n$, round sphere $\S^n$ or cylinder $\R\times \S^{n-1}$.
\end{thm}

\begin{rem}\label{pinching} 
Clearly, one cannot expect such results to hold without any assumption on the Weyl tensor $W$. For example, the classification of Einstein manifolds (such as $\mathbb{C}P^n$) with nonnegative scalar curvature, constant Q-curvature, and nonnegative Ricci curvature is a long-standing unsolved problem.

Actually, we can replace the assumption of vanishing Weyl tensor with one of the following pinching conditions:
	\begin{enumerate}
		\item[(1)] $R^2|W|\leq d_n|\nabla R|^2$ and Bach flat,
		\item[(2)]$R^2|W|+R|\nabla^2W|\leq \frac{2(n-3)}{2n-5}d_n|\nabla R|^2$,
	\end{enumerate}
	where $d_n= \frac{(n-2)(n^3-6n^2+16n-8)}{64n(n-1)^2}$.
	 Then $(M,g)$ is Ricci flat or isometric to  a quotient of round sphere $\S^n$ or cylinder $\R\times \S^{n-1}$, see Remark \ref{WandB}.
\end{rem}

Moreover, if we assume nonnegative sectional curvature instead of just nonnegative Ricci curvature, we can use a similar argument as in Theorem \ref{Thm 1.1} to establish the following rigidity theorem without requiring the locally conformally flat condition. To begin, we define the divergence of the Weyl tensor (see \eqref{curvature decomposition}) using normal coordinates,
\begin{align}\label{div weyl}
	(\delta W)_{ijk}=W_{ijkl,l}\qquad \mathrm{and}\qquad(\delta^2 W)_{ij}=W_{ikjl,lk}.
\end{align}

\begin{thm}\label{Thm 1.2}
Let $ (M^n,g)$ be a closed, oriented $n$-dimensional Riemannian manifold with $n\geq6$.  Suppose that 
$$\delta^2 W=0, \quad \quad sec\geq 0\quad\quad \mathrm{and}\quad\quad\int_{M}\nabla R\cdot\nabla \mathrm{Q}\ud V_g\leq0.$$ Then $(M^n,g)$ is   isometric to a quotient of the product of Einstein manifolds.
\end{thm}
\begin{rem}
It is not possible to directly derive Theorem \ref{Thm 1.2} from Theorem \ref{Thm 1.1}, as the latter does not satisfy the Weyl vanishing condition.
\end{rem}

Using Theorem \ref{Thm 1.2}, we can obtain some rigidity properties of compact immersed hypersurfaces in the standard model spaces $N^{n+1}(c)$. In this paper, we use $N^{n+1}(c)$ to denote $\S^{n+1}$, $\R^{n+1}$, and $\H^{n+1}$ when $c=1,0,-1$, respectively. The study of hypersurfaces in $N^{n+1}(c)$ has a long history. In 1977, Cheng and Yau \cite{S. Cheng and S.T. Yau} investigated hypersurfaces with constant scalar curvature in $N^{n+1}(c)$, and in particular, they proved the following results.

\begin{thm B}[\protect{Cheng and Yau \cite{S. Cheng and S.T. Yau}}]\label{yau thm}
If $(M,g)$ is a closed hypersurface in $N^{n+1}(c)$ with constant scalar curvature  $R\geq n(n-1)c$
	and  non-negative sectional curvature, 
then $M$ is either a totally umbilical hypersurface,  a product of two totally umbilical constant curved submanifolds or a flat manifold.
\end{thm B}

Motivated by Cheng and Yau's rigidity results, we obtain the following corollary from Theorem \ref{Thm 1.2}: 
\begin{cor}\label{surface}
	Assuming that $(M^n,g)$ is a closed and oriented hypersurface in $N^{n+1}(c)$ with $n\geq6$, and satisfies $\delta^2 W=0$, $sec\geq0$, and $\int_{M}\nabla R\cdot\nabla \mathrm{Q}\ud V_g\leq0.$
	\begin{enumerate}
	\item[(1)] If $c=0$, then $(M,g)$  is isometric to round sphere $\S^n$ that is  totally umbilical.
	\item[(2)] If $c=1$, then $(M,g)$ is isometric to either $\S^n$ or $\S^k\times\S^{n-k}$, $1\leq k \leq n-1$, and it is either totally umbilical or a Riemannian product of two totally umbilical submanifolds with constant curvature.
\item[(3)] If $c=-1$, then $(M,g)$  is isometric to metric sphere  that is  totally umbilical.
	\end{enumerate}
\end{cor}

Moreover, if we consider some extrinsic conditions instead of the non-negative sectional curvature assumption, motivated by \cite{H.Li}, we obtain the following result:
\begin{cor}\label{surface thm2}
Let $(M^n,g)$ be a compact, oriented hypersurface in $N^{n+1}(c)$, where $c=0$ or $1$, with $n\geq 6$. Assuming that $\delta^2W=0$ and $\int_{M}\nabla R\cdot\nabla \mathrm{Q}\ud V_g\leq0$, if the second fundamental form satisfies $\frac{H^2}{n}\leq|h|^2\leq\frac{H^2}{n-1}$, then $(M,g)$ must be isometric to a round sphere that is totally umbilical.
\end{cor}

Analogously to the Yamabe problem, a natural question is whether Q-curvature can be deformed to be  constant by   conformal change, which corresponds to solving a fourth-order elliptic partial differential equation  on the conformal factor,
\begin{equation}\label{l}
P_gu=cu^{\frac{n+4}{n-4}},\quad\quad u>0,
\end{equation}
where $P_g$ is referred to as the Paneitz operator, i.e,
\begin{equation}\label{Paneitz}
P_gu={\Delta}^2u-\mathrm{div}\left(\frac{((n-2)^2+4)R}{2(n-1)(n-2)}g(\nabla  u,\cdot)-\frac{4}{n-2} \mathrm{Ric}(\nabla u,\cdot)\right)+\frac{n-4}{2}\mathrm{Q}u.
\end{equation}
The equation \eqref{l} has been extensively studied, and numerous results have been obtained in \cite{M. Gursky F. Hang and Y.J. Lin, M. Gursky and A. Malchiodi, F.Hang and P. Yang 1, F.Hang and P. Yang 2, E. Hebey and F.Robert, Y.Y. Li and J. Xiong, J. Qing and D. Raske}.

Considering the analogy between the constant   Q-curvature problem and the Yamabe problem, despite substantial progress regarding the existence of solutions, to our knowledge, the issue of uniqueness (or lack thereof) has only been inspected in a small number of geometric settings. Actually, there are infinitely many branches of metrics with constant Q-curvature but without constant scalar curvature, as shown in \cite{BPS21}. Therefore, a natural question arises: is the metric with both constant scalar curvature and constant Q-curvature in a fixed conformal class unique up to scaling? In particular, we have the following uniqueness theorem. To be clearer, we replace $R$ with $R_g$ to distinguish between scalar curvatures with different metrics. Similarly, $Q_g$ is also metric-dependent.
\begin{thm}\label{Thm 1.3}
Let $(M^n,g)$ be a closed Riemannian manifold that is not conformally equivalent to the standard round sphere. Assuming that $R_g$ and $\mathrm{Q}_g$ are both constant, if $\tilde{g}$ is a metric in the conformal class $[g]$  such that $R_{\tilde{g}}$ and $\mathrm{Q}_{\tilde{g}}$ are also constant, then $\tilde{g}=c^2g$ for some positive constant $c$.
\end{thm}

It is worth pointing out that if the initial metric is only assumed to have constant scalar curvature (which can always be assumed by the Yamabe problem), the above result may not hold. To illustrate this, we provide the following example.

\begin{ex}
	Let $M = \S^1(T)\times \S^{n-1}$ and $L(\S^1(T))=T>2\pi(n-2)^{-1/2}$ be equipped with the standard product metric $\tilde{g}$. Then, there exists a metric $g=u^{\frac{4}{n-2}}\tilde{g}$ with a non-constant function $u$ such that
	$$R_{g}\equiv C_1>0,\quad \mathrm{Q}_{g}\not\equiv \mathrm{const};\quad \quad R_{\tilde{g}}\equiv C_2>0, \quad \mathrm{Q}_{\tilde{g}}\equiv \mathrm{const}. $$
\end{ex}

Finally, we would like to point out a rigidity result concerning metrics with constant scalar curvature and constant Q-curvature in the locally conformally flat case.

\begin{thm}\label{thm 4}
Let $(M,g)$ be a  connected locally conformally flat manifold with constant scalar curvature and constant Q-curvature, if there exist $p\in M$ such that $\mathrm{Ric}(p)>0$. Then  $(M^n,g)$ is isometric to a quotient of round sphere $\S^n$.
\end{thm}

The organization of this paper is as follows. In Section \ref{sec2}, we establish a formula for the Laplacian of the Ricci norm squared, which is important for the proof of Theorem \ref{Thm 1.1} and Theorem \ref{Thm 1.2}. In Section \ref{sec 3}, we provide a complete proof of Theorem \ref{Thm 1.1} and Theorem \ref{Thm 1.2}. In Section \ref{sec 4}, we use Theorem \ref{Thm 1.2} to prove Corollary \ref{surface} and Corollary \ref{surface thm2}. In Section \ref{sec 5}, we use a result by R. Schoen to complete the proof of Theorem \ref{Thm 1.3}. Finally, we prove Theorem \ref{thm 4}.

	\section{Preliminaries}\label{sec2}
We will begin this section by introducing some notation in Riemannian geometry using local coordinates. First, we select normal coordinates $(x_1,\cdots,x_n)$ centered at $p$, and define $\partial_i=\frac{\partial}{\partial x_i}$. Then, for any $\partial_i, \partial_j, \partial_k,$ and $\partial_l\in T_pM$, we define the Riemann curvature tensor as
\begin{align*}
R(\partial_i,\partial_j)\partial_k=\nabla_{\partial_j}\nabla_{\partial_i} \partial_k-\nabla_{\partial_i}\nabla_{\partial_j} \partial_k+\nabla_{[\partial_i, \partial_j]}\partial_k,
\end{align*}
and $R_{ijkl}=R(\partial_i,\partial_j,\partial_k,\partial_l)=\langle R(\partial_i,\partial_j)\partial_k,\partial_l\rangle_g$.

Using the contraction of the second and fourth indices of the Riemann curvature tensor, we can define the Ricci curvature as
\begin{align*}
R_{ij}={\rm Ric}(\pa_i,\pa_j)=\sum_{k,l=1}^{n}g^{kl}R_{ikjl},
\end{align*}
and contracting the Ricci curvature tensor once more, we obtain the scalar curvature $R=\sum_{i,j=1}^{n}R_{ij}g^{ij}$. We then define the trace-free part of the Ricci curvature as $$\overset{\circ}{\mathrm{Ric}}=\mathrm{Ric}-\frac{R}{n}g.$$
The decomposition of the curvature tensor is a well-known result and can be expressed as follows:
\begin{align}\label{curvature decomposition}
R_{ijkl}=A_{ij}\odot g_{kl}+W_{ijkl}
:=A_{ik}g_{jl}+A_{jl}g_{ik}-A_{il}g_{jk}-A_{jk}g_{il}+W_{ijkl}.
\end{align}
Here, the Schouten tensor is denoted by $A_{ij}$ and is defined as
 $$A_{ij}=\frac{1}{n-2}\left(R_{ij}-\frac{R_g}{2(n-1)}g_{ij}\right).$$  The Weyl tensor is denoted by $W_{ijkl}$ and is defined by \eqref{curvature decomposition}.

Another important tensor in Riemannian geometry is the Bach tensor, denoted by $B_{ij}$, which is defined as
\begin{equation}\label{Bach-Definition-1}
B_{ij}=\frac{1}{n-3}(\delta^2W)_{ij}+\frac{1}{n-2}R_{kl}W_{ikjl},
\end{equation}
where $\delta^2W$ is defined by \eqref{div weyl}.

To establish our main results, we begin by considering the quantity $|{\rm Ric}|^2$, which is the squared norm of the Ricci curvature tensor, and its Laplacian. Our motivation for this approach stems from two sources. Firstly, $|{\rm Ric}|^2$ appears in the Q-curvature equation \eqref{Q-curvature}. Secondly, if $|{\rm Ric}|^2$ is constant, it implies that ${\rm Ric}$ is parallel under certain geometric assumptions, and this, in turn, implies the desired geometric rigidity. This approach serves as a key ingredient in the proof of our main theorem.
	\begin{lem}\label{formula}
		Suppose $(M^n,g)$ is a complete  Riemannian manifold and $n>3$ , then
		\begin{equation}\label{Laplace-Ricci-Bach-Formula-1}
		\begin{split}
		\Delta |\mathrm{Ric}|^2
		=&2R_{ij,k}^2+\frac{n-2}{(n-1)}R_{ij}R_{,ij}+\frac{1}{(n-1)}R\Delta R\\[6pt]
		&+2I-4W_{ikjl}R_{kl}R_{ij}+2(n-2)B_{ij}R_{ij},
		\end{split}
		\end{equation}
		where  
		\begin{align}\label{I}
			I=\frac{n}{n-2}R_{ik}R_{jk}R_{ij}-\frac{2n-1}{(n-2)(n-1)}R|\mathrm{Ric}|^2+\frac{1}{(n-1)(n-2)}R^3.
		\end{align}
In particular,
\begin{enumerate}
\item [(1)] if $(M^n,g)$ is locally conformally flat, i.e., $W=0,$ then we have 
 \begin{equation}\label{Lap-Ricci-Zer-Wey-For-3}
			\Delta |\mathrm{Ric}|^2
			=2R_{ij,k}^2+\frac{n-2}{(n-1)}R_{ij}R_{,ij}+\frac{1}{(n-1)}R\Delta R+2I.
	 \end{equation}
\item [(2)]if $\delta^2W=0$, then we have 
 \begin{equation}\label{Lap-Ricci-Div-Wey-For-1}
 \begin{split}
 	\Delta |\mathrm{Ric}|^2
 = &2R_{ij,k}^2+\frac{n-2}{(n-1)}R_{ij}R_{,ij}+  \frac{1}{(n-1)}R\Delta R \\[6pt]
 &  +2R_{ij}R_{jk}R_{ki}-	2R_{ijkl}R_{ik}R_{jl}.
 \end{split}
	 \end{equation}
\end{enumerate}
\end{lem}
	\begin{pf}
By direct computation, we get
	\begin{equation} \label{Laplace-Ricci-Square-1}
		\begin{split}
		\Delta|\mathrm{Ric}|^2
		=&2R_{ij}\Delta R_{ij}+2R_{ij,k}^2.
		\end{split}
		\end{equation}
In the following, we will calculate the Laplacian of the Ricci tensor. This formula is well-known and can be found in the literature, for example, in \cite{CGY02} or \cite{Tian} in the case of dimension four. Here, we will provide a direct calculation for the convenience of the reader.

The main ingredients are the Bianchi
identities. To begin with, we have 
\begin{equation}
R_{ikjl,l}=R_{ij,k}-R_{jk,i},
\end{equation}
and 
\begin{equation}\label{Schouten-Div-For-1}
A_{ik,k}=\frac{1}{n-2}\left(R_{ik,k}-\frac{R_{,i}}{2(n-1)}\right)=\frac{1}{2(n-1)}R_{,i}.
\end{equation}

With the curvature  decomposition \eqref{curvature decomposition}, by the second  Bianchi identity again, 
we have 
\begin{equation}\label{Div-Weyl-For-1}
\begin{split}
W_{ikjl,l}&=
R_{ikjl,l}-\big(A_{ij,l}g_{kl}+A_{kl,l}g_{ij}-A_{il,l}g_{jk}-A_{jk,l}g_{il}\big)\\
&=(n-3)(A_{ij,k}-A_{jk,i}).
\end{split}
\end{equation}
Taking divergence of both sides of  \eqref{Div-Weyl-For-1}, we obtain
\begin{equation}\label{2Div-Weyl-For-1}
\begin{split}
(\delta^2W)_{ij}&=W_{ikjl,lk}=(n-3)(A_{ij,kk}-A_{jk,ik}).
\end{split}
\end{equation}
We recall a standard formulas for commuting
covariant derivatives, 
\begin{equation}\label{Schouten-Ric-For-1}
A_{jk,ik}=A_{jk,ki}-(R_{ikjl}A_{lk}-R_{il}A_{jl}).
\end{equation}
Therefore, combining \eqref{Schouten-Div-For-1}, \eqref{2Div-Weyl-For-1}  and \eqref{Schouten-Ric-For-1},  we have 
\begin{equation}\label{Laplace-Ric-2Div-W-1}
\begin{split}
\Delta A_{ij}=A_{ij,kk}&=\frac{1}{n-3}(\delta^2W)_{ij}+\frac{1}{2(n-1)}R_{,ij}-R_{ikjl}A_{lk}+R_{il}A_{jl}.
\end{split}
\end{equation}
In terms of  Bach tensor  \eqref{Bach-Definition-1},  \eqref{Laplace-Ric-2Div-W-1} can be expressed as  
\begin{equation}\label{Laplace-Schouten-Tensor-1}
\Delta A_{ij}=B_{ij}-\frac{1}{n-2}W_{ikjl}R_{kl}+\frac{1}{2(n-1)}R_{,ij}-R_{ikjl}A_{kl}+R_{il}A_{jl}.
\end{equation}
By the curvature  decomposition \eqref{curvature decomposition}, 
we get 
\begin{equation}\label{Cur-Dec-Cor-1}
 \begin{split}R_{ikjl}R_{kl}
&=\frac{1}{n-2}\Big(\frac{nR}{n-1}(R_{ij}-\frac{R}{n}g_{ij})+|\mathrm{Ric}|^2g_{ij}-2R_{ik}R_{kj}\Big)+W_{ikjl}R_{kl}.
\end{split}
\end{equation}
Therefore, we  may rewrite the   equation \eqref{Laplace-Schouten-Tensor-1} as
\begin{equation}\label{Laplace-Ricci-Tensor-3}
\begin{split}
\Delta R_{ij}&=(n-2)B_{ij}-2W_{ikjl}R_{kl}+\frac{n-2}{2(n-1)}R_{,ij}+\frac{1}{2(n-1)}\Delta Rg_{ij}\\
\quad &+\frac{1}{n-2}\left(\frac{R^2}{n-1}g_{ij}-\frac{nR}{n-1}R_{ij}-|\mathrm{Ric}|^2g_{ij}+nR_{ik}R_{kj}\right).\\[6pt]
\end{split}\end{equation}
Plugging  \eqref{Laplace-Ricci-Tensor-3} into \eqref{Laplace-Ricci-Square-1}, we get \eqref{Laplace-Ricci-Bach-Formula-1}.

In particular, if $(M^n,g)$ is locally conformally flat, i.e., $W=0$,  then \eqref{Lap-Ricci-Zer-Wey-For-3} is direct consequence of \eqref{Laplace-Ricci-Bach-Formula-1}.

While if $\delta^2W=0$, from \eqref{Bach-Definition-1},  then	\begin{equation}\label{Bach-DivWeyl-1}
    B_{ij}=\frac{1}{n-2}R_{kl}W_{ikjl}.
\end{equation}
Moreover, by   \eqref{Cur-Dec-Cor-1}, we have 
	\begin{align}\label{lig}
		W_{ikjl}R_{ij}R_{kl}
		=&R_{ikjl}R_{ij}R_{kl}-R_{ik}R_{kj}R_{ji}+I.
	\end{align}
Plugging  \eqref{Bach-DivWeyl-1} and \eqref{lig} into  \eqref{Laplace-Ricci-Bach-Formula-1}, we arrive  \eqref{Lap-Ricci-Div-Wey-For-1}.

\end{pf}

With the equations provided in Lemma \ref{formula}, proving Theorem \ref{Thm 1.1} hinges on a meticulous analysis of the non-negativity of the term $I$ \eqref{I}. To achieve this, we must first establish an algebraic inequality \eqref{Fundamental-inequality-2}. In this context, we present an alternative proof that enables us to identify the equality case using our method.

\begin{lem}\label{fundamental inequality}
	If $x_i\geq0,~\sum_{i=1}^{n}x_i=1$ and $n\geq 3$, then

\begin{equation}\label{Fundamental-inequality-2}n\sum_{i=1}^{n}x_i^3+\frac{1}{n-1}\geq\frac{2n-1}{n-1}\sum_{i=1}^{n}x_i^2.
\end{equation}
	The equality holds if and only if one of the following conditions is met:
	$(1)$ All $x_i$ are equal, i.e., $x_i=\frac{1}{n}$ for $i=1,\dots,n$.
	$(2)$ One of the $x_i$ is zero, and the remaining $n-1$ $x_i$ are equal, i.e., $x_i=\frac{1}{n-1}$ for $i=1,\dots,n-1$, and $x_n=0$ (up to permutation).
\end{lem}
\begin{pf}
We argue it by induction on $n$. 
	For $n=3$, we
	let $x_3=1-\sum_{i=1}^{2}x_i$ and 
	\begin{align*}
	    f_3(x_1,x_{2})&=3\left(\sum_{i=1}^{2}x_i^3+\left(1-\sum_{i=1}^{2}x_i\right)^3\right)-\frac{5}{2} \left(\sum_{i=1}^{2}x_i^2+\left(1-\sum_{i=1}^{2}x_i\right)^2\right)+\frac{1}{2}.
	\end{align*}
	If $x_i>0$ for all $i=1,2,3$,
	then all critical points of $f_3$ are given by $(\frac{1}{3},\frac{1}{3})$,$(\frac{4}{9},\frac{4}{9})$,$(\frac{4}{9},\frac{1}{9})$,$(\frac{1}{9},\frac{4}{9})$, and the corresponding critical values are 
	$$f_3\left(\frac{1}{3},\frac{1}{3}\right)=0,\quad f_3\left(\frac{4}{9},\frac{4}{9}\right)>0,\quad 
	f_3\left(\frac{4}{9},\frac{1}{9}\right)>0,
	\quad f_3\left(\frac{1}{9},\frac{4}{9}\right)>0.$$
If there exists an index $i$ such that $x_i=0$, without loss of generality, we can assume $x_3=0$ by symmetry. In this case, 
	\[f_3=f_3(x_1,1-x_{1})=4\left(x_1-\frac{1}{2}\right)^2\geq 0. \]
The equality holds if and  only if $x_1=\frac{1}{2}$, $x_2=\frac{1}{2}$, $x_3=0$.	 
	
If $n\geq4$,	we can assume that the conclusion holds for  $n-1$. Let
\begin{align*}
   	&f_n(x_1,\dots,x_{n-1})\\
   	=&n\left(\sum_{i=1}^{n-1}x_i^3+\left(1-\sum_{i=1}^{n-1}x_i\right)^3\right)-\frac{2n-1}{n-1} \left(\sum_{i=1}^{n-1}x_i^2+\left(1-\sum_{i=1}^{n-1}x_i\right)^2\right)+\frac{1}{n-1},
\end{align*}
	and the critical point $x^0=(x_1^0,x_2^0,\cdots,x_{n-1}^0)$ of $f_n$ satisfies $\frac{\partial f}{\partial x_i}(x^0)=0,i=1,\dots,n-1,$ i.e, 
	\[
	\left(x_i^0-\alpha\right)     \left(\left(x_i^0+\alpha\right)-\beta_n\right)=0,
	\]
	where $\alpha=1-\sum_{j=1}^{n-1}x_j^0,\beta_n=\frac{2(2n-1)}{3n(n-1)}$.
	The critical point also mets
	\begin{align*}
		\sum_{i=1}^{n-1}(x_i^0)^2=&\beta_n(1-\alpha)+\left(\alpha^2-\alpha\beta_n\right)(n-1),\\
		\sum_{i=1}^{n-1}(x_i^0)^3=&\beta_n\left(\beta_n(1-\alpha)+\left(\alpha^2-\alpha\beta_n\right)(n-1)\right)+\left(\alpha^2-\alpha\beta_n\right)(1-\alpha).
	\end{align*}
	Hence, we get
	\begin{align}\label{gd}
		f_n\left(x_1^0,\dots,x_{n-1}^0\right)=&n\left[\beta_n^2(1-\alpha)+\left(\alpha^2-\alpha\beta_n\right)\beta_n(n-1)+\left(\alpha^2-\alpha\beta_n\right)(1-\alpha)+\alpha^3\right]\no\\
		+&\frac{1}{n-1}-\frac{2n-1}{n-1}\left[\beta_n(1-\alpha)+\left(\alpha^2-\alpha\beta_n\right)(n-1)+\alpha^2\right]\no\\
		=&\left(n-\frac{n^2}{2}\beta_n\right)\left(\alpha^2-\beta_n\alpha\right)+\frac{1}{n-1}-\frac{n}{2}\beta_n^2.
	\end{align}
	Without loss of generality, we assume $$x_i^0=\alpha,\quad i=1,\dots,m;\quad\quad\quad x_i^0=\beta_n-\alpha,\quad i=m+1,\dots,n-1.$$
	By the definition of $\alpha$, we get
	\begin{align}\label{pkl}
		1-\alpha=m\alpha+(n-1-m)(\beta_n-\alpha),
	\end{align}
	then
	\begin{align}\label{dgt}
		\alpha=&\frac{\beta_n}{2}+\frac{n-2}{3(n-1)(2m-(n-2))}.
	\end{align}
We observe that if $2m=n-2$, then from \eqref{pkl}, we would have $\frac{n}{2}\beta_n=1$, which would further imply $n=2$. However, this contradicts the assumption that $n\geq3$. Therefore, by \eqref{dgt}, we can conclude that the minimum value of $f(x_1^0,\dots,x_{n-1}^0)$ is attained when $m=0$ or $m=n-1$. Additionally, we can infer that
	\begin{align*}
		m&=n-1,\quad\quad\quad x_i^0=\frac{1}{n},\quad i=1\dots,n-1;\\
		m&=0,\quad\quad\quad x_i^0=\frac{1}{n-1}-\frac{1}{3n(n-1)},\quad i=1,\dots,n-1.
	\end{align*}
	 Using \eqref{gd} and \eqref{dgt}, we know $$f_n\left(\frac{1}{n},\dots,\frac{1}{n}\right)=0,\quad f_n\left(\frac{1}{n-1}-\frac{1}{3n(n-1)},\dots,\frac{1}{n-1}-\frac{1}{3n(n-1)}\right)>0.$$ This indicates that any possible negative minimum must occur on the boundary. Without loss of generality, we can assume that $x_n=0$, and it suffices to check
	\begin{align}\label{ikj}
		n\sum_{i=1}^{n-1}x_i^3+\frac{1}{n-1}\geq\frac{2n-1}{n-1}\sum_{i=1}^{n-1}x_i^2.
	\end{align}
	By induction, we know
	\begin{align}\label{sgd}
		n\sum_{i=1}^{n-1}x_i^3+\frac{n}{(n-1)(n-2)}\geq&\frac{n(2n-3)}{(n-1)(n-2)}\sum_{i=1}^{n-1}x_i^2.
	\end{align}
	Now, plugging \eqref{sgd} into \eqref{ikj} we calculate
	\begin{align*}
		n\sum_{i=1}^{n-1}x_i^3+\frac{1}{n-1}=&n\sum_{i=1}^{n-1}x_i^3+\frac{n}{(n-1)(n-2)}-\frac{2}{(n-1)(n-2)}\\
		\geq&\frac{n(2n-3)}{(n-1)(n-2)}\sum_{i=1}^{n-1}x_i^2-\frac{2}{(n-1)(n-2)}\\
		=&\frac{2n-1}{n-1}\sum_{i=1}^{n-1}x_i^2+\frac{2}{n-2}\sum_{i=1}^{n-1}x_i^2-\frac{2}{(n-1)(n-2)}\\
		\geq&\frac{2n-1}{n-1}\sum_{i=1}^{n-1}x_i^2.
	\end{align*}
	The last step follows from the H\"older inequality. The equality holds if and only if $x_i=\frac{1}{n-1},i=1,\dots,n-1$, and $x_n=0$.
\end{pf}

Using Lemma \ref{fundamental inequality}, we can draw the following conclusion from \eqref{Lap-Ricci-Zer-Wey-For-3}.

\begin{lem}\label{key inequality}
Suppose $(M^n,g)$ is a complete and locally conformally flat Riemannian manifold. If $n>3$ and $\mathrm{Ric}\geq0$, then 
	\begin{equation}\label{Lap-Ricci-Zero-Weyl-PR-1}
		\Delta |\mathrm{Ric}|^2\geq2|\nabla \mathrm{Ric}|^2+
		\frac{n-2}{(n-1)}\langle \mathrm{Ric},\nabla^2R\rangle+\frac{1}{(n-1)}R\Delta R.
	\end{equation}
	
\end{lem} 
\begin{pf} 
	Given $p\in M$, we can choose normal coordinates centered at $p$ such that $R_{ij}(p)=\lambda_i(p)\delta_{ij}$. At $p$, the term $I$ in \eqref{Lap-Ricci-Zer-Wey-For-3} can be expressed as follows
\[
		I=\frac{1}{n-2}\left(n\sum_{i=1}^{n}\lambda_i^3-\frac{2n-1}{(n-1)}\left(\sum_{i=1}^{n}\lambda_i\right)\left(\sum_{i=1}^{n}\lambda_i^2\right)
		+\frac{1}{(n-1)}\left(\sum_{i=1}^{n}\lambda_i\right)^3\right).
\]
If $R(p)=0$, then $\mathrm{Ric}\geq 0$ implies that $\lambda_i(p)=0$ for all $i$, which in turn implies that $I=0$.
If $R(p)>0$, then by Lemma \ref{fundamental inequality} with $\lambda_i(p)\geq 0$ for all $i$, we have $I\geq 0$. Therefore, \eqref{Lap-Ricci-Zero-Weyl-PR-1} follows immediately from \eqref{Lap-Ricci-Zer-Wey-For-3} and $I\geq 0$.
\end{pf}

\section[Proof of Theorem \ref{Thm 1.1} and Theorem 1.2]{Two rigidity theorems for Riemannian  manifolds}\label{sec 3}

In this section, using the Lemma \ref{formula}, Lemma \ref{fundamental inequality} and Lemma \ref{key inequality}, we prove the Theorem \ref{Thm 1.1} and Theorem \ref{Thm 1.2}.\\

\textbf{Proof of Theorem \ref{Thm 1.1}}: 
To simplify the notation, we define $a_n=\frac{1}{2(n-1)}$, $b_n=\frac{2}{(n-2)^2}$, and $c_n=\frac{n^3-4n^2+16n-16}{8(n-1)^2(n-2)^2}$. Using the definition of Q-curvature  \eqref{Q-curvature}, we can derive that
		\begin{align*}
		-\Delta \mathrm{Q}=a_n\Delta^2R-c_n\Delta R^2+b_n\Delta|\mathrm{Ric}|^2.
		\end{align*}
Since $M$ is a Riemannian manifold that is locally conformally flat and satisfies $\mathrm{Ric}\geq0$, we can combine this with Lemma \ref{key inequality} to obtain the following result,
		\begin{align*}
		-\Delta \mathrm{Q}-a_n\Delta ^2R+c_n\Delta R^2\geq& 2b_nR_{ij,k}^2+
		\frac{n-2}{(n-1)}b_nR_{ij}R_{,ij}+\frac{b_n}{(n-1)}R\Delta R.
		\end{align*}
	Multiplying both sides of the inequality by $R$ and integrating over $M$ yields:
\begin{equation}\label{1.3}
	\begin{split}
	\int_{M}\nabla R\cdot&\nabla \mathrm{Q}\ud V_g-a_n\int_{M}(\Delta R)^2\ud V_g-2c_n\int_{M}|\nabla R|^2R\ud V_g\\
	&\geq 2b_n\int_{M}|\nabla\mathrm{Ric}|^2R\ud V_g+\frac{n-2}{(n-1)}b_n\int_{M}R_{ij}R_{,ij}R\ud V_g\\
	&-\frac{2b_n}{n-1}\int_{M}|\nabla R|^2R\ud V_g.
	\end{split}
	\end{equation}
For the second term on the right hand side of \eqref{1.3}, integrating by parts, we have 
	\begin{align}\label{is}
	\int_{M}R_{ij}R_{,ij}R\ud V_g&=	\int_{M}(R_{ij}R_{,i}R)_{,j}-R_{ij,j}R_{,i}R-R_{ij}R_{,i}R_{,j}\ud V_g\no\\
	&=-\frac{1}{2}\int_{M}|\nabla R|^2R\ud V_g-\int_{M}\mathrm{Ric}(\nabla R,\nabla R)\ud V_g.
	\end{align}
Recall the Bochner formula
\[\frac{1}{2}\Delta|\nabla R|^2=|\nabla^2R|^2+\langle \nabla R, \nabla\Delta R\rangle_g+\mathrm{Ric}(\nabla R,\nabla R),\]
	we  have 
	\begin{align*}
	\int_{M}\mathrm{Ric}(\nabla R,\nabla R)+|\nabla^2R|^2\ud V_g=\int_{M}(\Delta R)^2\ud V_g,
	\end{align*}
	therefore
	\begin{align}\label{oj}
	\int_{M}(\Delta R)^2\ud V_g\geq\frac{n}{n-1}\int_{M}\mathrm{Ric}(\nabla R,\nabla R)\ud V_g.
	\end{align}
	Moreover, it is easy to see that
	\begin{equation}\label{Gra-Ric-Gra-Sca-Ine-1}
	|\nabla \mathrm{Ric}|^2\geq \frac{1}{n}|\nabla R|^2.
	\end{equation}
	Plugging \eqref{is}, \eqref{oj} and \eqref{Gra-Ric-Gra-Sca-Ine-1} into \eqref{1.3}, we have
%
	\begin{align}\label{sljd}
	\int_{M}\nabla R\cdot\nabla \mathrm{Q}\ud V_g&+	\left(\left(\frac{n-2}{2(n-1)}+\frac{2}{n-1}-\frac{2}{n}\right)b_n-2c_n\right)\int_{M}|\nabla R|^2R\ud V_g\no\\
	&+\left(\frac{n-2}{n-1}b_n-\frac{n}{n-1}a_n\right)\int_{M}\mathrm{Ric}(\nabla R,\nabla R)\ud V_g\geq 0.
	\end{align}
Note that for any $n\geq 3$, one can get
	\begin{align}\label{sjd}
-l_n:=\left(\frac{n+2}{2(n-1)}-\frac{2}{n}\right)b_n-2c_n=-\frac{n^3-6n^2+16n-8}{4n(n-1)^2(n-2)}<0.
	\end{align}
While if  $n\geq6$, the  coefficient
	\begin{align}\label{Dimension-condition-1}
	\frac{n-2}{n-1}b_n-\frac{n}{n-1}a_n&=\frac{-n^2+6n-4}{2(n-1)^2(n-2)}<0.
	\end{align}
Furthermore, our assumption that $	\int_{M}\nabla R\cdot \nabla \mathrm{Q}\ud V_g\leq 0$, then \eqref{sljd} implies that
\begin{equation}\label{Scalar-Ricci-Parallel-1}
\int_{M}|\nabla R|^2R\ud V_g=0,\quad\quad\int_{M}\mathrm{Ric}(\nabla R,\nabla R)\ud V_g=0.
\end{equation}
Plugging \eqref{Scalar-Ricci-Parallel-1} into \eqref{1.3}, with \eqref{is},
we have  
\begin{equation}\label{1.4}
	0\geq 2b_n\int_{M}|\nabla \mathrm{Ric}|^2R\ud V_g+a_n\int_{M}(\Delta R)^2\ud V_g.
\end{equation}

Now we know that $R$ is a constant. There are two cases will happen.
\begin{enumerate}
\item[(1)]If $R\equiv0$, then $\mathrm{Ric}\equiv0$. Since $(M, g)$ is locally conformally flat, then it is flat, i.e, a quotient of $\R^n$.  
\item[(2)]If $R\equiv \mathrm{const}>0$, then using \eqref{1.4}, we obtain $\nabla \mathrm{Ric}\equiv0$, which means that the Ricci tensor is parallel. Therefore, by the well known de Rham-Wu's splitting theorem  \cite[Theorem 10.3.1]{Petersen}, it follows that $M$ is isometric, at least locally, to a product of Einstein manifolds.

To be more precise, let us assume that $M$ is simply connected (otherwise, we consider the universal covering of $M$) and decompose $T_pM=\oplus_i V_{i,p}$ into the orthogonal eigenspaces $V_{i,p}={\rm Ker}({\rm Ric}-\lambda_i(p) {\rm Id})$ for ${\rm Ric}$ with respect to distinct eigenvalues $\lambda_1<\cdots <\lambda_m$ at point $p$. Since ${\rm Ric}$ is parallel, we can parallel translate these eigenspaces to obtain a global decomposition of the tangent bundle $TM=V_1\oplus\cdots\oplus V_m$ into parallel distributions, with the property that ${\rm Ric}|_{V_i}=\lambda_i(p){\rm Id}$. In particular, the eigenvalues of ${\rm Ric}$ are constant on $M$.

By the de Rham-Wu's splitting theorem,  we know
	$$(M,g)=(M_1\times\cdots\times M_m, g_1+\cdots+g_m)$$
	with  $TM_i=V_i$.
	We assert that $(M_i,g_i)$ is an Einstein manifold for $i=1,\dots,n$. To simplify the argument, we will only consider the case where $i=1$. We select an orthonormal basis $\{e_\alpha\}_{\alpha=1}^{n} \subset TM$ such that $\{e_i\}_{i=1}^{m_1} \subset TM_1$. It is worth noting that $(M_1,g_1)$ is totally geodesic with respect to $(M,g)$. Therefore,
	\begin{align}\label{ric}
	\mathrm{Ric}_{g_1}(e_i, e_j)=&\sum_{k=1}^{m_1}Rm_{g_1}(e_i,e_k,e_j,e_k)=\sum_{k=1}^{m_1}Rm_{g}(e_i,e_k,e_j,e_k)\no\\
	=&\sum_{\alpha=1}^{n}Rm_{g}(e_i,e_{\alpha},e_j,e_{\alpha})=\mathrm{Ric}_{g}(e_i,e_j)=\lambda_1\delta_{ij},
	\end{align}
	
	where $Rm_g(\cdot,\cdot,\cdot,\cdot)$ denotes the Riemann curvature tensor of metric $g$, the last equality of \eqref{ric} follows by the the definition of $V_1$. So, we complete the claim. 
	
	Lemma \ref{fundamental inequality} implies that $m\leq2$. If $m=1$, then $(M,g)$ is a constant sectional curvature manifold, either $\R^n$ or $\S^n$. If $m=2$, then $dimV_1=n-1$ and $dimV_2=1$, so we can express $(M,g)$ as a product manifold $(M',g_1)\times(N,g_2)$, where $N$ is a one-dimensional flat manifold. Since $(M,g)$ is locally conformally flat, we claim that $(M',g_1)$ is also a locally conformally flat Einstein manifold. With the same notation as above $m_1=n-1$, for $1\leq i,j,k,l\leq n-1$
	\begin{align}\label{Weyl}
	W_{ijkl}^{'}=&R_{ijkl}^{'}-\frac{1}{n-3}\left(R_{ij}^{'}-\frac{R_g^{'}}{2(n-2)}g_{ij}^{'}\right)\odot g_{kl}^{'}\no\\
	=&R_{ijkl}-\frac{1}{n-3}\left(R_{ij}-\frac{
		R_{g}}{2(n-2)}g_{ij}\right)\odot g_{kl}\no\\
	=&R_{ijkl}-\frac{1}{n-2}\left(R_{ij}-\frac{R_{g}}{2(n-1)}g_{ij}\right)\odot g_{kl} -\frac{R_{ij}}{(n-3)(n-2)}\odot g_{kl}\no\\
	&\quad\quad+\frac{R_{g}}{(n-2)(n-3)(n-1)}g_{ij}\odot g_{kl}.
	\end{align}
	Due to $W_{ijkl}=0$ and $R_{ij}=\lambda_1 g_{ij}=\frac{R_g}{n-1}g_{ij}$, it follows from equation \eqref{Weyl} that $W_{ijkl}^{'}=0$. (We also note that the above argument only holds for $dim M^{'}=n-1$). Since we know that $(M,g)$ is simply connected, it follows that both $(N, g_2)$ and $(M',g_1)$ are also simply connected. Therefore, $(N,g_2)$ is isometric to $\R$ and $(M',g_1)$ is isometric to a positively curved constant sectional curvature manifold, namely $\S^{n-1}$.
	
	\end{enumerate}

	\begin{rem}\label{WandB}
If we drop the Weyl vanishing condition, the general formula \eqref{Laplace-Ricci-Bach-Formula-1} implies that an extra term involving $W$ will appear in \eqref{sljd}. Specifically, the `0' on the right-hand side of \eqref{sljd} will be replaced by
		\begin{align*}
		II: =2b_n\int_{M}((n-2)B_{ij}R_{ij}-2W_{imjk}R_{mk}R_{ij})R\ud V_g.
		\end{align*}
As mentioned in Remark \ref{pinching}, 	
	if  $(M,g)$ satisfies pinching condition  (1) or (2),   i.e., 		
		\begin{enumerate}
			\item[(1)] If $R^2|W|\leq d_n|\nabla R|^2$ and Bach flat, then
			\begin{align*}
			II\geq& -4b_n\int_{M}|W||\mathrm{Ric}|^2R\ud V_g\\
			\geq&-\frac{l_n}{2}\int_{M}|\nabla R|^2R\ud V_g.
			\end{align*}
			\item[(2)] If $R^2|W|+R|\nabla^2W|\leq \frac{2(n-3)}{2n-5}d_n|\nabla R|^2$, then
			\begin{align*}
			II\geq& -2b_n\int_{M}|W||\mathrm{Ric}|^2R\ud V_g-\frac{2(n-2)}{n-3}b_n \int_{M}|\nabla^2W||\mathrm{Ric}|R\ud V_g\\
			\geq& -2b_n\int_{M}|W|R^3\ud V_g-\frac{2(n-2)}{n-3}b_n \int_{M}|\nabla^2W|R^2\ud V_g\\
			\geq &-\frac{l_n}{2}\int_{M}|\nabla R|^2R\ud V_g,
			\end{align*}
		\end{enumerate}
	where $l_n$ is given by \eqref{sjd}. In either case, the above argument still holds, and one can conclude that $\nabla\mathrm{Ric}=0$. Moreover, if $R=0$, then $(M,g)$ is Ricci flat. If $R>0$, by applying pinching condition (1) or (2), we obtain $W\equiv0$. In particular, the same conclusion as in Theorem \ref{Thm 1.1} still holds.\\
	\end{rem}

\textbf{Proof of Theorem \ref{Thm 1.2}}: For any $p\in M$, we choose coordinates such that $R_{ij}(p)=\lambda_{i}\delta_{ij}$, then
	 \begin{align}\label{dk}
	 	R_{ij}R_{jk}R_{ki}-	R_{ijkl}R_{ik}R_{jl}=&R_{ii}\lambda_{i}^2-R_{ijij}\lambda_{i}\lambda_{j}=\frac{1}{2}R_{ijij}(\lambda_{i}-\lambda_{j})^2.
	 \end{align}
In the case where $\delta^2W=0$, using \eqref{dk}, we can rewrite equation \eqref{Lap-Ricci-Div-Wey-For-1} as
	 \begin{align}\label{Laplace-Ricci-
	Div-Weyl-1}
	 	\Delta|\mathrm{Ric}|^2
	 =&2R_{ij,k}^2+\frac{1}{(n-1)}R\Delta R+\frac{n-2}{(n-1)}R_{ij}R_{,ij}+R_{ijij}(\lambda_i-\lambda_j)^2.
	 \end{align} 
If the sectional curvature is nonnegative, i.e. $R_{ijij}\geq 0$,  from \eqref{Laplace-Ricci-
	Div-Weyl-1}, we have 
	 \begin{align}\label{imp}
	 	\Delta|\mathrm{Ric}|^2
	 	 \geq&2R_{ij,k}^2+\frac{1}{(n-1)}R\Delta R+\frac{n-2}{(n-1)}R_{ij}R_{,ij}.
	 \end{align} 
Furthermore, if $n\geq6$ and $\int_{M}\nabla R\cdot\nabla \mathrm{Q}\ud V_g\leq0$, then using \eqref{imp} and the proof of Theorem \ref{Thm 1.1}, we obtain $\nabla \mathrm{Ric}=0$ and $R=const$. Therefore, by the well-known de Rham-Wu's splitting theorem \cite[Theorem 10.3.1]{Petersen}, $M$ is isometric, at least locally, to a product of Einstein manifolds.

\section{Some rigidity results of hypersurface in model spaces}\label{sec 4}
In this section, we will explore some rigidity properties of compact immersed hypersurfaces in the standard model spaces $N^{n+1}(c)$, where  $N^{n+1}(c)$ to denote $\S^{n+1}$, $\R^{n+1}$, and $\H^{n+1}$ when $c=1,0,-1$, respectively. Our main focus is on proving that the principal curvatures are constant, which allows us to classify the hypersurface. \\

\textbf{Proof of the Corollary \ref{surface}}: 
Using Theorem \ref{Thm 1.2}, we can conclude that $\nabla \mathrm{Ric}=0$, i.e., the hypersurface $(M,g)$ has a parallel Ricci tensor and constant scalar curvature. Although we have obtained a lot of information about the intrinsic curvature, in order to classify the hypersurface, we need to know the principal curvatures. For any $p\in M$, we can choose coordinates such that the second fundamental form $h$ takes the form $h_{ij}=\kappa_i\delta_{ij}$. The Gauss formula tells us that
	 \begin{equation}\label{gauss equation}
	    R_{ijkl}=(\delta_{ik}\delta_{jl}-\delta_{il}\delta_{jk})c+h_{ik}h_{jl}-h_{il}h_{jk}=(c+\kappa_i\kappa_j)(\delta_{ik}\delta_{jl}-\delta_{il}\delta_{jk}),
	 \end{equation}
	 \begin{equation}\label{gauss ric}
	     R_{ij}=(n-1)c\delta_{ij}+Hh_{ij}-h_{ik}h_{jk}=((n-1)c+H\kappa_i-\kappa_i^2)\delta_{ij},
	 \end{equation}
	 \begin{equation*}\label{gauss scalar}
	      R=n(n-1)c+H^2-\|h\|_g^2=n(n-1)c+\left(\sum_{i=1}^{n}\kappa_i\right)^2-\sum_{i=1}^{n}\kappa_i^2.
	 \end{equation*}
	 
	\begin{enumerate}
	 \item[(1)] In $\R^{n+1}$, if $R=0$, then the embedding $(M,g) \hookrightarrow \R^{n+1}$ is flat. However, this leads to a contradiction as there must be a point $p_0$ on the compact hypersurface in $\R^{n+1}$ where $h(p_0)>0$. Therefore, we can assume $R>0$, which implies $H^2>0$. Next, we can assume $H>0$. We define $\lambda_i=H\kappa_i-\kappa_i^2\geq0$, which implies that $\kappa_i\geq0$. Combining equations \eqref{gauss equation} and \eqref{gauss ric}, we obtain
\begin{align}\label{dhy}
	 	R_{ijkl}R_{ik}R_{jl}-R_{ij}R_{jk}R_{ki}=&\sum_{i\not=j}\kappa_i\kappa_j\lambda_i\lambda_j-\sum_{i=1}^{n}\lambda_i^3\no\\
	 	=&\left(\sum_{i=1}^{n}\kappa_i\lambda_i\right)^2-\sum_{i=1}^{n}\kappa_i^2\lambda_i^2-\sum_{i=1}^{n}\lambda_i^3\no\\
	 	=&\left(\sum_{i=1}^{n}\kappa_i\lambda_i\right)^2-H\sum_{i=1}^{n}\kappa_i\lambda_i^2\leq 0.
	 \end{align}
	 So, combining \eqref{dk} and \eqref{Laplace-Ricci-
	 	Div-Weyl-1}, we can deduce that the above inequality \eqref{dhy} must hold as an equality. By the Cauchy inequality, for any $p$ in $M$, there exists a natural number $m_p$ such that $\lambda_1=\cdots=\lambda_{m_p}>0$ and $\kappa_{m_p+1}=\cdots\kappa_n=0$. From $\nabla \mathrm{Ric}=0$, we can infer that the eigenvalues of $\mathrm{Ric}$ are constant. Therefore, we have $m_p=m$ and $\lambda_1=\cdots=\lambda_{m}=\frac{R_g}{m}>0$, $\lambda_{m+1}=\cdots\lambda_n=0$, and $\kappa_i>0$ for $1\leq i\leq m$.
	 
	 If $m=n$, then we claim that $\kappa_1=\cdots=\kappa_n=\kappa$. Otherwise, if there exist distinct $\kappa_i$ and $\kappa_j$, then we can write
	 \begin{align}\label{wuss}
	 H\kappa_i-\kappa_i^2=H\kappa_j-\kappa_j^2,
	 \end{align}
	 then \eqref{wuss} implies that $H=\kappa_i+\kappa_j$. This yields that $\kappa_k=0$ for $k\not=i,j$, it is impossible  for $m=n$.
	 
	  In $\R^{n+1}$, we know that there exists a point $p_0$ such that $h(p_0)>0$.   If $m<n$, but $\kappa_{m+1}=\cdots=\kappa_n=0$. This leads to a contradiction. Therefore, $(M,g)$ is an umbilical hypersurface in $\R^{n+1}$ and must be a sphere.

	 \item[(2)]In $\S^{n+1}$, we define $\lambda_i=n-1+H\kappa_i-\kappa_i^2\geq0$, and from \eqref{Laplace-Ricci-
	 	Div-Weyl-1}, we know that if $\lambda_i\not=\lambda_j$, then $R_{ijij}=1+\kappa_i\kappa_j=0$. We claim that there are at most two distinct values of $\lambda_i$. If not, suppose there exist three distinct $\lambda_i$, $\lambda_j$, and $\lambda_k$. Then, $\kappa_i$, $\kappa_j$, and $\kappa_k$ must also be different from each other. But,
	 \begin{align}\label{sdhd}
	 \begin{cases}
	 \displaystyle 1+\kappa_i\kappa_j=&0,\\
	\displaystyle 1+\kappa_k\kappa_j=&0,\\
\displaystyle	1+\kappa_i\kappa_k=&0.
	 \end{cases}
	 \end{align}
	 Taking the difference of each equation\eqref{sdhd}, you will obtain $\kappa_i=\kappa_j$. Without loss of generality, we let $$\lambda_1=\cdots=\lambda_m=\lambda,\quad\quad\lambda_{m+1}=\cdots=\lambda_n=\mu.$$
	 The following statement asserts that $(M,g)$ is necessarily isoparametric with a maximum of 2 distinct principal curvatures, denoted as $\kappa_i$. This means that $\kappa_i$ must be constant for $i=1,2,\cdots,n$, and only two unique values of $\kappa_i$ are permitted.
	 \begin{enumerate}	
	 	 \item[(I)] If $m=0$ or $n$, then $(M,g)$ is Einstein, but we only consider the case where  $m=n$. For  any $p\in M$, there can be  at most two distinct values of principal curvature $\kappa_i$. If not,  there exists  three different values, say $\kappa_i$, $\kappa_j$, $\kappa_k$. Then,
	 \begin{align}\label{H and kappa}
	 \begin{cases}
	 \displaystyle n-1+H\kappa_i-\kappa_i^2=&\lambda\\
	  \displaystyle n-1+H\kappa_j-\kappa_j^2=&\lambda\\
	  \displaystyle n-1+H\kappa_k-\kappa_k^2=&\lambda
	 \end{cases}
	 \end{align}
	Taking the difference of each equation \eqref{H and kappa}, we get $H=\kappa_i+\kappa_j=\kappa_j+\kappa_k$, which contradicts the assumption of three distinct values of principal curvatures. If there exist $p\in M$ such that $\kappa_1=\cdots=\kappa_n$, then $\lambda_i\geq n-1+(n-1)\kappa_i^2\geq n-1$, using the Page 200 of Cheng and Yau \cite{S. Cheng and S.T. Yau} we are done. Thus, we can assume there exists $1\leq m_p\leq n-1$, and we also assume $\kappa_1=\cdots=\kappa_{m_p}=\kappa(p)$ and $\kappa_{m_p+1}=\cdots\kappa_n=t(p)$. Then, we know $H=\kappa(p)+t(p)=m_p\kappa(p)+(n-m_p)t(p)$ and $\kappa(p)t(p)=const$, where the second identity follows from the fact that $\lambda=n-1+(\kappa(p)+t(p))\kappa(p)-\kappa(p)^2$ is constant. Hence,
	 \begin{align}\label{k and t}
	 \begin{cases}
	 \displaystyle	(m_p-1)\kappa(p)+(n-m_p-1)t(p)=0\\
	 \displaystyle	\kappa(p)t(p)=const:=d.
	 	\end{cases}
	 \end{align}
	 If there exist $p\in M$ such that $\kappa(p)=0$ or $t(p)=0$, then  $\lambda=n-1$ and $R=n(n-1)$, by \cite{S. Cheng and S.T. Yau} we are done. Now, we assume that for any $p\in M$,  $\kappa(p)\not=0$ and $t(p)\not=0$, then $m_p\not=1$ and $n-1$. So, the above equations \eqref{k and t} turn into
	 $$(m_p-1)\kappa(p)^2+(n-m_p-1)d=0,$$
	where $\kappa(p)$ only takes discrete values and is a continuous function above $p$, which implies that $m_p=m$ and $\kappa(p)\equiv\kappa$, $t(p)\equiv t$ are both constant and satisfies
	 \begin{align*}
	 \begin{cases}
	 \displaystyle (m-1)\kappa+(n-m-1)t=0,\\
	 \displaystyle n-1+\kappa t=\lambda.
	 \end{cases}
	 \end{align*}
	In fact, $(M,g)$ is a isoparametric hypersurface, Cartan's formula indicates that $kt+1=0$. So,
	\begin{align*}
		\lambda=n-2,\quad\quad k=\pm\sqrt{\frac{n-m-1}{m-1}},\quad\quad t=\mp\sqrt{\frac{m-1}{n-m-1}}.
	\end{align*}
	 \item[(II)] If $1\leq m<n$,  we focus on $\kappa_i$. Firstly, we have $\kappa_i\not=\kappa_{l}$ for $1\leq i\leq m$ and $m+1\leq l\leq n$. Secondly, we claim that for any $p\in M$ $\kappa_1(p)=\cdots=\kappa_m(p)=\kappa(p)$ and $\kappa_{m+1}(p)=\cdots=\kappa_n(p)=t(p)$. If  there exists $\kappa_i(p)\not=\kappa_j(p)$ for $1\leq i,j\leq m$, then
	 \begin{align}\label{djfu}
	 \begin{cases}
	 \displaystyle 1+\kappa_i(p)\kappa_l(p)=&0,\\
	 \displaystyle 1+\kappa_j(p)\kappa_l(p)=&0.
	 \end{cases}
	 \end{align}
	 Hence, we get $\kappa_i(p)=\kappa_j(p)$ from \eqref{djfu} , which is a contradiction. Similarly, we can obtain $\kappa_{m+1}(p)=\cdots=\kappa_n(p)$. If $m\not=1$, then $1+\kappa(p)t(p)=0$ and $\lambda=n-1+(m\kappa(p)+(n-m)t(p))\kappa(p)-\kappa(p)^2$, we obtain $\kappa(p)$ is constant. If $m=1$, then $\lambda=n-1+(\kappa(p)+(n-1)t(p))t(p)-t(p)^2$, we obtain $t(p)$ is constant. In both case, $\kappa$ and $t$ are constants and satisfy $$1+\kappa t=0.$$
	 
	In conclusion, in both cases, we have $\nabla \mathrm{Ric}=0$, and $\kappa_i$ are constant. Moreover, there exists $m\in \N$ such that $\kappa_1=\cdots=\kappa_m=\kappa$ and $\kappa_{m+1}=\cdots=\kappa_n=t$. Then, by applying the classification of isoparametric hypersurfaces in $\S^{n+1}$ with at most two different principal curvatures \cite[Page 99]{Cecil Ryan}, we obtain that $(M,g)$ is isometric to $\S^k\times \S^{n-k}$.
	 	 \end{enumerate}
	 \item[(3)] In the hyperbolic space $\H^{n+1}$, we consider the quantities $\lambda_i = H\kappa_i - \kappa_i^2 - (n-1) \geq 0$. If $\lambda_i \neq \lambda_j$, then $R_{ijij} = \kappa_i \kappa_j - 1 = 0$. In this case, we can assume that $H > 0$ and $\kappa_i \geq 0$ for $1 \leq i \leq n$, and $H$ and $\kappa_i$ have the same sign. By a similar argument as before, we conclude that $(M,g)$ is an isoparametric hypersurface in $\H^{n+1}$.
	 
	 Therefore, by applying \cite[Theorem 3.14, page 97]{Cecil Ryan}, we obtain that $(M,g)$ is either a totally geodesic hyperplane, an equidistant hypersurface, a horosphere, a metric hypersphere, or a tube over a totally geodesic submanifold of codimension greater than one in $\H^{n+1}$. However, only the metric sphere is compact in hyperbolic space. Hence, $(M,g)$ must be a metric sphere.\\
	 \end{enumerate}

\textbf{Proof of the Corollary \ref{surface thm2}}: As stated in Corollary \ref{surface}, we define the second fundamental form as $h_{ij}=\kappa_i\delta_{ij}$, where $\kappa_i$ represents the principal curvatures.
Let $Z=h-\frac{H}{n}{\rm I}$ and $\mu_i=\kappa_i-\frac{H}{n}$. It follows that $\sum_{i}\mu_i=0$. For all $1\leq i \leq n$, the following holds,
\begin{equation}\label{eigenvalue-estimate-1}
|Z|^2=\sum_{j=1,j\neq i}^n\mu_j^2+\mu_i^2\geq\frac{1}{n-1}\left(\sum_{j=1,j\neq i}^n\mu_j\right)^2+\mu_i^2=\frac{n}{n-1}\mu_i^2.
\end{equation}
We now see that the hypothesis of  $\frac{H^2}{n}\leq|h|^2\leq\frac{H^2}{n-1}$ implies that
\begin{equation}\label{Pinching-eig-est-1}
|Z|^2=|h|^2-\frac{H^2}{n}\leq \frac{H^2}{n(n-1)}.
\end{equation}
Combining \eqref{eigenvalue-estimate-1} and 
\eqref{Pinching-eig-est-1},  for any $1\leq i\leq n$, we have
\begin{align}\label{estimate of mu}
 |\mu_i|\leq\frac{|H|}{n}.
\end{align}
Consequently, with \eqref{estimate of mu}, we have (pointwise) 
\begin{enumerate}
\item[(1)] either $H> 0$, $\kappa_i=\mu_i+\frac{H}{n}\geq 0$, for $1\leq i \leq n$,
\item[(2)] or $H<0$, $\kappa_i=\mu_i+\frac{H}{n}\leq 0$, for $1\leq i \leq n$,
\item[(3)] or $H=0$, $\kappa_i=\frac{H}{n}= 0$, for $1\leq i \leq n$.
\end{enumerate}
In any scenario, the inequality $R_{ijij}=c+\kappa_i\kappa_j\geq 0$ holds true for $c=0$ or $1$, respectively. By combining this with Corollary \ref{surface} and the pinching condition $\frac{H^2}{n}\leq|h|^2\leq\frac{H^2}{n-1}$, we can conclude that $M$ is isometric to a round sphere that is totally umbilical. It is worth noting that examples of Riemannian products of umbilical round spheres are ruled out by the pinching condition. Therefore, the proof of Corollary \ref{surface thm2} is now complete.

\section{Constant Q-curvature and constant scalar curvature}\label{sec 5}

Let $(M^n,g)$ be a closed Riemannian manifold with $n>2$. We will now prove the uniqueness, up to scaling, of metrics with constant scalar curvature and constant Q-curvature within a fixed conformal class.

Firstly, we recall the Obata Theorem:
\begin{thm}[\protect{M.Obata \cite{M.Obata}}]\label{Obata thm}
	If $(M,g)$ is Einstein and $\tilde{g} \in [g]$ has constant scalar curvature, then
	\begin{itemize}
		\item[(1)]  $\tilde{g}$ is Einstein.
		\item[(2)]Furthermore, if $(M,g)$ is not conformally equivalent to the standard round sphere, then there exists a positive constant $c$ such that $\tilde{g}=c^2g$.
	\end{itemize}	
\end{thm}  
In the general case, R. Schoen proved the following result, as stated in (1.12) of \cite{R.Schoen}. Note that all geometry notations will have a subscript $g$ added as a subscript to distinguish between different metrics.

\begin{lem}\label{conformal property}
	If $R_{g}\equiv C$, then $\forall\, \tilde{g} \in [g], ~\tilde{g}=e^{2f}g$ we have 
	
	\begin{align*}
		\int_{M}e^{-f}|\overset{\circ}{\mathrm{Ric}_g}|_g^2 \mathrm{d}V_g \leq
		 \int_{M} e^{-f} | \overset{\circ}{\mathrm{Ric}_{\tilde{g}}} |_{g}^2 \mathrm{d}V_g.
	\end{align*}
	
\end{lem}
\begin{pf}
	By the conformal change of Ricci tensor, we obtain
	\begin{align}\label{con1}
	\widetilde{\mathrm{Ric}}_{ij}=&\mathrm{Ric}_{ij}-(n-2)f_if_j+(n-2)f_{ij}-(\Delta_gf+(n-2)|\nabla_gf|_g^2)g_{ij}
	\end{align}
	and
	\begin{align}\label{con2}
	\tilde{R}_{\tilde{g}}=&e^{-2f}\left(R_g-2(n-1)\Delta_gf-(n-1)(n-2)|\nabla_gf|_g^2\right).
	\end{align}
	Thus, using \eqref{con1} and \eqref{con2} we have
	\begin{align*}
	\overset{\circ}{\widetilde{\mathrm{Ric}}}_{ij}
	=&\widetilde{\mathrm{Ric}}_{ij}-\frac{\tilde{R}_{\tilde{g}}}{n}\tilde{g}_{ij}\\
	=&\mathrm{Ric}_{ij}-(n-2)f_if_j+(n-2)f_{ij}-(\Delta_gf+(n-2)|\nabla_gf|_g^2)g_{ij}\\
	&\quad\quad~-\frac{1}{n}(R_g-2(n-1)\Delta_gf-(n-1)(n-2)|\nabla_gf|_g^2)g_{ij}\\
	=&\overset{\circ}{\mathrm{Ric}_{ij}}-(n-2)f_if_j+(n-2)f_{,ij}+\frac{n-2}{n}\Delta_gfg_{ij}-\frac{n-2}{n}|\nabla_gf|_g^2g_{ij}.
	\end{align*}
	This yields
	\begin{align}\label{con3}
	&e^{-f}\langle\overset{\circ}{\widetilde{\mathrm{Ric}}},\overset{\circ}{\mathrm{Ric}}\rangle_g\no\\
	=&e^{-f}\left(\overset{\circ}{\mathrm{Ric}_{ij}}-(n-2)f_if_j+(n-2)f_{ij}+\frac{n-2}{n}\Delta_gfg_{ij}-\frac{n-2}{n}|\nabla_gf|_g^2g_{ij}\right)\overset{\circ}{\mathrm{Ric}_{kl}}g^{ik}g^{jl}\no\\
	=&e^{-f}\left(|\overset{\circ}{\mathrm{Ric}_g}|_g^2-(n-2)\overset{\circ}{\mathrm{Ric}}^{ij}f_if_j+(n-2)\overset{\circ}{\mathrm{Ric}}^{ij}f_{,ij}\right)\no\\
	=&e^{-f}|\overset{\circ}{\mathrm{Ric}_g}|_g^2+(n-2)\left(e^{-f}\overset{\circ}{\mathrm{Ric}}^{ij}f_i\right)_{,j};
	\end{align}
	the last equality follows from 
	\begin{align*}
	\overset{\circ}{\mathrm{Ric}}^{ij},j={\mathrm{Ric}}^{ij},j-\frac{(R_g)_{,j}}{n}g^{ij}=\frac{n-2}{2n}(R_g)_{,j}g^{ij}=0.
	\end{align*}
	Combining with \eqref{con3}, we obtain
	\begin{align*}
	\int_{M}e^{-f}|\overset{\circ}{\mathrm{Ric}_g}|_g^2\mathrm{d}V_g=&\int_{M}e^{-f}\langle\overset{\circ}{\widetilde{\mathrm{Ric}}},\overset{\circ}{\mathrm{Ric}}\rangle_g\mathrm{d}V_g\\
	\leq &\left(\int_{M}e^{-f}|\overset{\circ}{\mathrm{Ric}_g}|_g^2\mathrm{d}V_g\right)^{1/2}\left(\int_{M} e^{-f} | \overset{\circ}{\mathrm{Ric}_{\tilde{g}}} |_{g}^2 \mathrm{d}V_g\right)^{1/2}.
	\end{align*}
If $\overset{\circ}{\mathrm{Ric}_g}\equiv 0$, then the Obata Theorem implies that $\overset{\circ}{\mathrm{Ric}_{\tilde{g}}}\equiv 0$, so the inequality is trivial.
If $\overset{\circ}{\mathrm{Ric}_g}\not\equiv 0$, then $\overset{\circ}{\mathrm{Ric}_{\tilde{g}}}\not\equiv 0$ as well, and the above inequality still supports our argument.
\end{pf}
Next, we will prove the uniqueness of metrics with constant scalar curvature and constant Q-curvature within a conformal class.\\

\textbf{Proof of Theorem \ref{Thm 1.3}}: 
	Since we know
\begin{align*}
	| \overset{\circ}{\mathrm{Ric}_{\tilde{g}}} |_{g}^2=e^{4f}| \overset{\circ}{\mathrm{Ric}_{\tilde{g}}} |_{\tilde{g}}^2 \qquad \mathrm{and}\qquad \mathrm{d}V_g=e^{-nf}\ud V_{\tilde{g}}, 
\end{align*}
we can rewrite	the Lemma \ref{conformal property} as
	\begin{align}\label{thm 3: formula a}
		\int_{M}e^{-f}|\overset{\circ}{\mathrm{Ric}_g}|_g^2 \mathrm{d}V_g \leq&
		\int_{M} e^{-(n-3)f} | \overset{\circ}{\mathrm{Ric}_{\tilde{g}}} |_{\tilde{g}}^2 \mathrm{d}V_{\tilde{g}}.
	\end{align}
	Due to the relation $g=e^{-2f}\tilde{g}$, replacing $f$ by $-f$, we can similarly obtain
	\begin{align}\label{thm 3: formula b}
		\int_{M}e^{f}| \overset{\circ}{\mathrm{Ric}_{\tilde{g}}} |_{\tilde{g}}^2 \mathrm{d}V_{\tilde{g}}\leq&
		\int_{M} e^{(n-3)f} |\overset{\circ}{\mathrm{Ric}_g}|_g^2 \mathrm{d}V_g.
	\end{align}
	By the definition of Q-curvature \eqref{Q-curvature}, we
 notice that the condition $$\mathrm{Q}_g\equiv \mathrm{const},~R_g\equiv \mathrm{const} $$
 is equivalent to 
 \begin{align*}
 	|\overset{\circ}{\mathrm{Ric}_g}|_g^2\equiv \mathrm{const},~R_g\equiv \mathrm{const}.
 \end{align*}
 Set $|\overset{\circ}{\mathrm{Ric}_g}|_g^2\equiv C_1$ and $| \overset{\circ}{\mathrm{Ric}_{\tilde{g}}} |_{\tilde{g}}^2\equiv C_2$, there are two cases will happen.

  \textbf{Case 1}: If $C_1=0$,  since $(M, g)$ is not conformally equivalent to the standard round sphere, we get $C_2=0$ and $\tilde{g}=c^2g$ based on the Obata Theorem \ref{Obata thm} .

 \textbf{Case 2}: 
 	 If $C_1\neq0$, we can conclude, according to Obata Theorem \ref{Obata thm}, that $C_2\neq0$. We would like to emphasize that in this case, the assumption of ``not conformally equivalent to the standard sphere" is unnecessary. Instead, we will directly apply integral and H\"older inequality to prove it. From \eqref{thm 3: formula a} and \eqref{thm 3: formula b}, we obtain
 	\begin{align}\label{con4}
 	C_1\int_{M}e^{-f}\mathrm{d}V_g \leq C_2\int_{M} e^{3f}\mathrm{d}V_g 
 	\end{align}
 	and
 	\begin{align}\label{con4.5}
 	C_2\int_{M}e^{(n+1)f}\mathrm{d}V_g\leq C_1\int_{M} e^{(n-3)f}\mathrm{d}V_g.
 	\end{align}
 	By applying	H\"older inequality and \eqref{con4.5}, we can obtain the following inequality
 	\begin{align}\label{con5}
 	\int_{M}e^{(n+1)f}\mathrm{d}V_g\leq&\left(\frac{C_1}{C_2}\right)^{(n+1)/4}\mathrm{Vol}_g(M),
 	\end{align}	
 	then \eqref{con4} and \eqref{con5} implies that
 	\begin{align}\label{dkwiw}
 	\mathrm{Vol}_g(M)^2
 	\leq&\int_{M}e^{-f}\mathrm{d}V_g\int_{M}e^{f}\mathrm{d}V_g\leq \left(\frac{C_2}{C_1}\right)\int_{M} e^{3f}\mathrm{d}V_g\int_{M}e^{f}\mathrm{d}V_g\no\\
 	\leq&\left(\frac{C_2}{C_1}\right)\left(\int_{M}e^{(n+1)f}\mathrm{d}V_g\right)^{\frac{3}{(n+1)}}\left(\mathrm{Vol}_g(M)\right)^{\frac{(n-2)}{(n+1)}}\no\\
 	&\quad\quad\quad\left(\int_{M}e^{(n+1)f}\mathrm{d}V_g\right)^{\frac{1}{(n+1)}}\left(\mathrm{Vol}_g(M)\right)^{\frac{n}{(n+1)}}\no\\
 	\leq&\left(\frac{C_2}{C_1}\right) \left(\left(\frac{C_1}{C_2}\right)^{\frac{(n+1)}{4}}\mathrm{Vol}_g(M)\right)^{\frac{3}{(n+1)}}\left(\mathrm{Vol}_g(M)\right)^{\frac{(n-2)}{(n+1)}}\no\\
 	&\quad\quad\quad\left(\left(\frac{C_1}{C_2}\right)^{\frac{(n+1)}{4}}\mathrm{Vol}_g(M)\right)^{\frac{1}{(n+1)}}\left(\mathrm{Vol}_g(M)\right)^{\frac{n}{(n+1)}}\no\\
 	\leq&\mathrm{Vol}_g(M)^2.
 	\end{align}	
 	Therefore, from \eqref{dkwiw}, we have
 	\begin{align*}
 	\mathrm{Vol}_g(M)^2
 	=\int_{M}e^{-f}\mathrm{d}V_g\int_{M}e^{f}\mathrm{d}V_g,
 	\end{align*}
 	we can use H\"older inequality together with the facts that $e^{-f}>0$ and $e^{f}>0$ to conclude that $e^f\equiv c$, where $c$ is a positive constant.

Let us consider the following problem: given an initial metric, we can apply the Yamabe problem to assume that it has constant scalar curvature. Then, if we make a proper conformal change such that the new metric has constant scalar curvature and constant $Q$-curvature, is the conformal factor also constant? If the scalar curvature of the initial metric is a nonpositive constant, the uniqueness of the Yamabe problem implies that the conformal factor is also a constant.

However, if the initial metric has a positive constant scalar curvature, there exists a counterexample. As previously mentioned, there are metrics with constant $Q$-curvature but without constant scalar curvature, as shown in \cite{BPS21}. Specifically, consider the initial metric $g$, and use it to construct a new metric that has constant scalar curvature but without constant Q-curvature.
\begin{ex}
	Let $M = \S^1(T)\times \S^{n-1}$ and $L(\S^1(T))=T>2\pi(n-2)^{-1/2}$ be equipped with the standard product metric $\tilde{g}$. Then, there exists a metric $g=u^{\frac{4}{n-2}}\tilde{g}$ with a non-constant function $u$ such that
	$$R_{g}\equiv C_1>0,\quad \mathrm{Q}_{g}\not\equiv \mathrm{const};\quad \quad R_{\tilde{g}}\equiv C_2>0, \quad \mathrm{Q}_{\tilde{g}}\equiv \mathrm{const} .$$
\end{ex}
\begin{pf}
The universal cover of $M$ is $\R\times \S^{n-1}$. Let $e_1\in T\R$ and $e_2,\dots, e_n\in T\S^{n-1}$ with $|e_i|=1$ for $i=1,\dots,n$. Then, we have $K_{\tilde{g}}(e_1,e_i)=0$ for $i\neq 1$, and $K_{\tilde{g}}(e_i,e_j)=1$ for $i\neq j$ and $i,j\geq2$. From these, we obtain
\begin{align*}
	\mathrm{Ric}_{\tilde{g}}(e_1,e_1)=&\sum_{i=2}^{n}K_{\tilde{g}}(e_1,e_i)=0,\quad\quad  \mathrm{Ric}_{\tilde{g}}(e_i,e_i)=\sum_{j\not=i,j=2}^{n}K_{\tilde{g}}(e_i,e_j)=n-2,~i\geq2;\\
\mathrm{Ric}_{\tilde{g}}(e_1,e_i)=&0,~i\geq2,\qquad\qquad\qquad\quad \mathrm{Ric}_{\tilde{g}}(e_i,e_j)=0,i\not=j,~i,j\geq2.
\end{align*}
Thus,
$$R_{\tilde{g}}=(n-1)(n-2),\quad\quad\quad \mathrm{Q}_{\tilde{g}}=\frac{n^3-4n^2}{8}.$$
We consider the  Yamabe equation:
$$-\frac{4(n-1)}{(n-2)}\Delta_{\tilde{g}}u+(n-1)(n-2)u=n(n-1)u^{\frac{n+2}{n-2}},\quad\quad \mathrm{\quad on\quad} M.$$
 If $T>2\pi(n-2)^{-1/2}$, a non-constant solution $u$ exists, as shown in \cite{R.Schoen}. By defining $g=u^{\frac{4}{n-2}}\tilde{g}$, we obtain $R_g=n(n-1)$. Since the initial metric $\tilde{g}$ is not Einstein, we can apply \textbf{Case 2}  from Theorem \ref{Thm 1.3} to establish that if $\mathrm{Q}_{g}\equiv \text{const}$, then $u\equiv \text{const}$. However, this conclusion leads to a contradiction.
\end{pf}

\textbf{Proof of Theorem \ref{thm 4}}: We define $U=\{p\in M\mid \mathrm{Ric}(p)>0\}$. The hypothesis of $\mathrm{Ric}(p)>0$ at some point $p$ ensures that $U$ is nonempty. Moreover, $U$ is obviously an open set.

Using \eqref{Lap-Ricci-Zer-Wey-For-3} and the assumption of Theorem \ref{thm 4}, we obtain
\begin{align*}
0=R_{ij,k}^2+I.
\end{align*}
For any $p\in U$, applying Lemma \ref{fundamental inequality}, we conclude that $I(p)=0$. Consequently, the eigenvalue of $\mathrm{Ric}$ must be constant, and $\nabla \mathrm{Ric}=0$ in $U$. In this case, we have

$$R_{ijkl}(p)=c(g_{ik}g_{jl}-g_{il}g_{jk})(p)$$ 
where $c>0$ is independent of $p$.

If $p_0\in \partial U$, then there exists a sequence $\{p_i\}_{i=1}^{\infty}\subset U$ such that $p_i\to p_0$. Using the continuity of the curvature tensor, we obtain

$$R_{ijkl}(p_0)=c(g_{ik}g_{jl}-g_{il}g_{jk})(p_0).$$
Therefore, $U$ is also a closed set. Since $(M,g)$ is a connected manifold, we conclude that $U=M$ and it has constant curvature.\\\\
	{\noindent\small{\bf Acknowledgment:} We deeply appreciate the referee for providing us with invaluable feedback and insightful suggestions that have significantly improved the quality of this paper.

{\noindent\small{\bf Data availability:} Data sharing not applicable to this article as no datasets were generated or analysed during the current study.
	\section*{Declarations}
	{\noindent\small{\bf Conflict of interest:} On behalf of all authors, the corresponding author states that there is no conflict of interest.

\bibliographystyle{unsrt}

	\bigskip
\noindent Y.Y. Xu

\noindent Department of Mathematics, Nanjing University, \\
Nanjing 210093, China\\[1mm]
Email: \textsf{xuyiyan@nju.edu.cn}\\

\noindent S. Zhang

\noindent Department of Mathematics, Nanjing University, \\
Nanjing 210093, China\\[1mm]
Email: \textsf{dg21210019@smail.nju.edu.cn}

\medskip  	

	\end{document}